\documentclass{conm-p-l}

\newtheorem{theorem}{Theorem}[section]
\newtheorem{lemma}[theorem]{Lemma}
\newtheorem{proposition}[theorem]{Proposition}

\theoremstyle{definition}

\theoremstyle{remark}
\newtheorem{remark}[theorem]{Remark}

\newcommand{\bbC}{{\mathbb C}}

\newcommand{\BbC}{{\mathbb C}}

\newcommand{\BbZ}{{\mathbb Z}}
\newcommand{\BbP}{{\mathbb P}}

\newcommand{\param}{H^2\setminus\{0\}}
\newcommand{\ma}{{\mathcal A}}
\newcommand{\mc}{{\mathcal C}}
\newcommand{\rone}{{\rm I}^{(n)}}
\newcommand{\rii}{{\rm II}^{(n)}}
\newcommand{\riii}{{\rm III}^{(n)}}
\newcommand{\riv}{{\rm IV}^{(n)}}
\newcommand{\scl}{(-\sqrt{-1}\hbar)}

\numberwithin{equation}{section}

\begin{document}

\title{Gauge Theory Techniques in Quantum Cohomology}

\author{Steven Rosenberg}
\address{Department of Mathematics and Statistics, Boston University, Boston, 
Mass\-a\-chu\-setts 02215}
\email{sr@math.bu.edu}

\author{Mihaela Vajiac}
\address{Department of Mathematics, Univeristy of Texas,
Austin, Texas}
\email{mbvajiac@math.utexas.edu}
\thanks{The second author was partially supported by the Clay Mathematics Institute Liftoff Program.}

\date{\today}


\begin{abstract}
Quantum cohomology   gives a finite dimensional
integrable system via the Dubrovin connection.  Motivated by
Givental's work on mirror symmetry, we
use gauge theory techniques and the
Frobenius Integrability Theorem to find flat sections for the
Dubrovin connection.  An explicit calculation is given for projective
space. 
\end{abstract}

\maketitle

\section{Introduction}

The work of Givental \cite{G}, \cite{P} and Liu-Lian-Yau \cite{LLY}
 on mirror symmetry 
relates Gromov-Witten
invariants of a quintic hypersurface in ${\mathbb P}^3$ to period
integrals of K\"ahler structures on the 
mirror manifold.  Givental's method uses detailed
calculations of equivariant GW invariants to produce flat sections of
the Dubrovin connection on the tangent bundle to the even cohomology
of the hypersurface,
which are then related to solutions of the Picard-Fuchs ODE for the
 periods on the mirror.

Givental's approach leads to the following general
question:  given a flat connection on the tangent bundle to a
vector space, how can we compute the flat sections?  
In contrast to  Donaldson/Seiberg-Witten/Chern-Simons gauge theories,
where gauge directions are quotiented out,
the moduli space of flat connections here is trivial. Thus all the
information in the Dubrovin connection is
 contained in the gauge
transformation taking the Dubrovin connection to the trivial connection.
Since the Dubrovin connection is defined in terms of the quantum
product, which in turn encodes GW
invariants, this gauge transformation captures all the algebra of the
quantum product.

The purpose of
this paper is to compute this gauge transformation and the
corresponding flat sections
from a classical PDE point of
view, specifically through a systematic use of the Frobenius
Integrability Theorem.  It is surprisingly difficult to
compute flat sections; even for the cup product, where all quantum
corrections are turned off, the Dubrovin connection is rather trivial,
but we show in \S2 that the flat sections are much more complicated.
In \S3, we give a general expression for the flat sections, which
involves exponentiating an infinite matrix whose entries are
matrices.  In \S4, we reduce the calculation of the flat sections for the
quantum product on ${\mathbb P}^m$ to exponentials of ordinary
matrices, and compare our formula with known
solutions.

In this paper, as in \cite{G}, we restrict attention to the small quantum
product, with some comments at the end on the big quantum product and
on coupling to gravity.  This last product, with its relation to
infinite dimensional integrable systems, is the real case of interest,
and the current paper was intended as a warm-up.  The unexpected 
complexity of the calculations in this case indicates the difficulty
of proving the higher genus Virasoro conjecture \cite{LT}.

We would like to thank David Fried for helpful conversations.

\section{Constant Products on Vector Spaces}

Let $\ast$ be  an associative, commutative product 
on a complex vector space ${\mathcal H}$.  
The associated  Dubrovin
connection on $T{\mathcal H}$ is
\begin{equation}\label{dub}\nabla_Y X=dX(Y)+\sqrt{-1}Y\ast X
\end{equation}
with connection one-form $\omega (Y)(X)=\sqrt{-1}Y\ast X$. If $\{
T_0,\ldots,T_m\}$ is a basis of ${\mathcal H}$ with 
\begin{equation}\label{multtable}T_i\ast T_j=\Gamma_{ij}^k
T_k,\end{equation}
 then $\omega_j^i =\sqrt{-1} \Gamma_{\ell j}^i\ dt^\ell$ (where a typical
$\alpha\in {\mathcal H}$ is $\alpha =t_iT^i$).  It is fundamental that
$\nabla$ is flat because the product is associative and commutative.
Note that $\nabla$ stays flat if we replace 
 $\sqrt{-1}$ in (\ref{dub}) by $\hbar\in{\BbC}$.

The space ${\mathfrak M}$ 
of flat connections modulo gauge transformations over any
manifold ${\mathcal H}$ has formal tangent space $H^1({\mathcal H}, {\mathfrak
gl}(m,{\BbC}))$ at $\nabla$, where dim ${\mathcal H} =m$ and the de Rham
cohomology is computed with respect to the exterior derivative coupled
to $\nabla$ (see e.g.~\cite[\S7.2]{A}).  For our ${\mathcal H}$, the cohomology vanishes by a
Poincar\'e Lemma, and the
points of ${\mathfrak M}$ are characterized by the holonomy of the connection.
Since ${\mathcal H}$ is contractible, there is no nontrivial holonomy, so
${\mathfrak M}$ consists of one point.  

In particular,
$\nabla$ is globally gauge
equivalent to the trivial connection $d$, when $d_Y(X)=dX(Y)$. Thus
there exists $g\colon \bbC^m\rightarrow {\rm
GL}(m,\bbC )$ with $g\cdot d\equiv g^{-1}dg=\nabla$, which translates
into $g^{-1}dg=\omega$. A basis of the flat sections for $\nabla$ is
given by $\{ g^{-1} e_i\}_{i=0}^m$, where $e_i=\partial_{t^i}$ is the
constant section, since $\nabla g^{-1}e_i=(g^{-1}dg)g^{-1}e_i=0$. This
shows

\begin{lemma}\label{one}
 Let $g$ be the gauge transformation on $T{\mathcal H}$ such
that $g^{-1}dg =\nabla$. Then a basis of the flat sections of $\nabla$
is given by the columns of $g^{-1}$.\end{lemma}

Thus finding flat sections is equivalent to solving $g^{-1}dg=\omega$
for $g$. In coordinates, this becomes the system
\begin{equation}
\frac{\partial g_j^i}{\partial t^\ell}=\sqrt{-1} g_s^i
\Gamma_{j \ell}^s,
\label{e:one}
\end{equation}
If $A\in {\rm GL}(m^2,\bbC )$ is considered as a constant
gauge transformation, then $g^{-1}dg=\omega $ iff
$(Ag)^{-1}d(Ag)=\omega$. Thus we may assume
$g(0,0,\cdots,0)={\rm Id}$ as an initial condition, which is natural
as $g^{-1}$ will then take flat sections $\{ e_i\}$ of $d$ at $t\in
T_0{\mathcal H}$ to flat sections of $\nabla$ through $t$.

(\ref{e:one}) is a classical PDE handled by the Frobenius
Integrability Theorem (see e.g.~\cite[Vol.~I, Ch.~6]{S}).  We can rewrite
(\ref{e:one}) as
\begin{equation}
\frac{\partial g_j^i}{\partial t^\ell}=f_\ell (t, g(t)),
\label{e:two}
\end{equation}
with $x,f_\ell (t,x)\ m\times m$ matrices satisfying
\begin{eqnarray*}
f_\ell (t,x)_j^i &=& \sqrt{-1} x_s^i \Gamma_{\ell j}^s\cr
&=& \sqrt{-1}(x\cdot \Gamma_\ell )_j^i ,
\end{eqnarray*}
where $\Gamma_\ell$ is the $m\times m$ matrix $(\Gamma_\ell )_j^i
=\Gamma_{\ell j}^i$.

For later purposes, we first
treat the case where $\ast$ is independent of
$t\in {\mathcal H}$. This includes the case of the classical cup
product. Here the
Frobenius integrability conditions for the system  (\ref{e:two}) are
\begin{equation}
\Gamma_{ik}^j \Gamma_{\ell m}^k =\Gamma_{\ell k}^j \Gamma_{im}^k .
\label{e:three}
\end{equation}
Since $\ast$ is commutative, $\Gamma_{jk}^i =\Gamma_{kj}^i$, and
(\ref{e:three} )
is equivalent to the associativity of $\ast$. In addition, (\ref{e:three}) is
equivalent to $[\Gamma_i,\Gamma_\ell ]=0$, which implies that the
exponentials of the $\Gamma_i$ commute (see below).

The proof of the Frobenius
Integrability Theorem leads to a construction of the solution of
(\ref{e:one}).  We
first solve the $m\times m$ matrix valued, complex time  ODE
\begin{eqnarray*}\dot B(t,0,\ldots,0) &=& f_1((t,0,\ldots,0),B(t,0,
\ldots, 0))=\sqrt{-1}B(t,0,\ldots,0)\cdot\Gamma_0,\\
B(0) &=& {\rm Id}.\end{eqnarray*}
Thus $B(t,0,\ldots, 0)=\exp [\sqrt{-1} t\Gamma_0]$. We then fix $t^0$
and solve
\begin{eqnarray*}\dot B(t^0,t,0,\ldots,0) &=&
f_1(t^0,t,0,\ldots,0,B(t^0,t,0, \ldots,0))
=\sqrt{-1}B(t)\cdot \Gamma_1,\\
B(0)&=&\exp [\sqrt{-1}t^0 \Gamma_0],\end{eqnarray*}
so $B(t^0,t,0,\ldots,0) =\exp [\sqrt{-1}t^0\Gamma_0]\exp
[\sqrt{-1}t^1 \Gamma_1]$. 
Continuing, we get at the last step
\[g(t^0,t^1,\ldots, t^m)
= \prod_{\ell =0}^m \exp [\sqrt{-1}t^\ell
\Gamma_\ell ]
=\exp [\sqrt{-1} \sum_{\ell =0}^m t^\ell \Gamma_\ell ]
\]
since the $\Gamma_\ell$  commute. This shows:

\begin{proposition} Let $\ast$ be a commutative, associative product
on a complex vector space ${\mathcal H}$.  Define matrices $\Gamma_i,
\ i =
1,\ldots,{\rm dim}\ {\mathcal H}$ by $T_i\ast T_j = (\Gamma_i)^k_j T_k$
for a basis $\{T_k\}$ of ${\mathcal H}$.  Then
a basis of the flat sections 
of the associated flat Dubrovin connection on $T{\mathcal H}$  is
$$\left\{ \exp [- \sqrt{-1}t^\ell \Gamma_\ell ]\right\}_{\ell =0}^m.
$$\end{proposition}

For example, if $\ast$ is the cup product on ${\mathbb P}^1$, then the
connection one-form on $T{\mathcal H} = TH^*({\mathbb P}^1,{\mathbb C})$ is
$$\omega =\sqrt{-1} \left( \begin{array}{cc}
dt^0&0\\
dt^1&dt^0\end{array}\right),\ \Gamma_0={\rm Id},\ \Gamma_1 =\left(
 \begin{array}{cc}0&0\\ 1&0 \end{array}\right),
$$
and the flat sections are spanned by the columns of
$$\exp  \left[ \sqrt{-1} \left( \begin{array}{cc}
t^0&0\\ t^1&t^0\end{array}\right)
\right] =\left( \begin{array}{cc}e^{-it^0}&0\\
-it^1 e^{-it^0} &e^{-it^0}\end{array} \right).
$$

\section{The Small Quantum Product}

   In this section, we will find  a basis for the flat sections
for the Dubrovin connection for the small quantum product. As
explained in e.g.~\cite{BCPP}, this is a deformation of the cup
product on the even cohomology
 ${\mathcal H} = H^{\rm ev}(M;{\mathbb C})$ of 
a symplectic manifold or smooth projective variety $M$, 
but with  deformations only in $H^2(M;{\mathbb C})$ directions.  We
restrict attention to weakly monotonic symplectic manifolds (e.g.~Fano
varieties), so that
there are no convergence issues for the quantum product.

More precisely, the small quantum product is defined   on $i^\ast
T{\mathcal H}$, where $i: H^2\rightarrow {\mathcal H}$ is the inclusion. For
$\{ T_1,\ldots,T_k\}$  a basis of $H^2$, the small quantum
product is defined as in (\ref{multtable}) by
$\Gamma_{jr}^i (t^1,\ldots,t^k)=\sqrt{-1}\psi_{jr\ell}
h^{\ell i}$, where $(h_{\alpha\beta})=(T_\alpha \cdot T_\beta )$ is
the intersection matrix, $(h^{\alpha\beta})$ its inverse, and
\begin{equation}\label{inserttwo}
\psi_{jr\ell} =T_j\cdot T_r\cdot T_\ell +\sum_{\beta}\exp
\left[ \sum_{s=1}^k t^s\int_\beta T_s\right]I_\beta (T_j,T_r,T_\ell
)\end{equation} 
\cite[\S2.1]{BCPP}.
In (\ref{inserttwo}),  $I_\beta (T_j,T_r,T_\ell )$ is the three-pointed
Gromov-Witten invariant, and $\beta$ ranges over $H_2(M,{\mathbb
Z})\setminus \{0\}.$

Since the product now depends on $(t^1,\ldots,t^k)\in H^2$, the
flatness of $\nabla$ is more complicated.  More precisely, we extend
(\ref{dub}) to a family of connections 
\begin{equation}\label{extdub}(\nabla_\hbar)_YX = dX(Y) + \hbar Y\ast
X,\end{equation} 
$\hbar\in \mathbb C$, with curvature $\hbar d\omega
+\hbar^2\omega\wedge\omega$, $\omega^i_j = \Gamma^i_{\ell j}dt^\ell.$
Then $\omega\wedge\omega =0$ because $\ast$ is commutative and
associative, and $d\omega =0$ because the small quantum product is
potential (i.e. $\omega$ is exact).  For convenience we will set
$\hbar = \sqrt{-1}$ in this section.

  
To find a basis of flat sections, we apply the Frobenius
  Integrability Theorem method.  We must first
solve the ODE for $t^1$:
  $$\frac{\partial g_j^i}{\partial t} =\sqrt{-1} g_r^i
[T_j\cdot T_1\cdot T_\ell +\sum_{\beta\in \param}\exp [t\int_\beta T_1]
I_\beta (T_j,T_1,T_\ell )]h^{\ell r}.
$$
This is of the form

\begin{equation}
\dot B(t)=B(t) \left[ A+\sum_{i=1}^n e^{a_it} C_i\right] 
\label{e:five_a}
\end{equation}
with $A_j^r =\sqrt{-1}(T_j\cdot T_1\cdot T_\ell )h^{\ell r}$,
$a_i=\int_{\beta_i} T_1$, and $(C_i)_j^r =I_{\beta_i}(T_j,T_1,T_\ell
)h^{\ell r}$. Here $\{ \beta_1,\ldots,\beta_n\}$ is the set of classes
in $H_2$ such that $I_\beta (T_j,T_1,T_\ell )\not= 0$ for some
$j,\ell$.  We may assume that $a_i\in {\mathbb Z}.$ 

Let $B^{(1)}(t)$ be the solution of (\ref{e:five_a}) with $B^{(1)}(0) =
{\rm Id}.$  As in \S2, we
will then solve similar equations for $t^2,\ldots, t^k$ and get
$$g(t) =
B^{(1)}(t^1,0,\ldots,0)B^{(2)}(t^1,t^2,0,\ldots,0)\cdot\ldots\cdot
B^{(k)}(t^1,\ldots,t^k).$$ 
We actually want the columns of $g^{-1}$, so we must invert 
$B =B^{(i)}.$  Since $\dot B^{-1} = -B^{-1}\dot B B^{-1}$, we have 
\begin{equation}\label{e:five}
 \dot B^{-1}(t) = -\left[ A+\sum_{i=1}^n e^{a_it}
C_i\right] B^{-1}(t).\end{equation}

Equations of the form (\ref{e:five_a}), (\ref{e:five}) (which differ
 only by taking the transpose)
 have a long history.  They arise in the
Riemann-Hilbert problem of constructing second order ODEs with regular
singular points having
prescribed monodromy at prescribed points in the complex plane.
Classically, (\ref{e:five}) was treated by a good ansatz leading to
complicated recursion formulas for the coefficients.  The main
technical contribution  of this paper is a more modern organization of
the recursion relations.

The rest of this section is devoted to solving (\ref{e:five})
with the assumption $n=1$ and $a_1=1$. For convenience, we replace
$A,C$ in (\ref{e:five}) 
by $-A, -C$, respectively. The general case is discussed
in \S5.

Under the substitution $x=e^t$, (\ref{e:five}) becomes
$\frac{dB}{dx} =-[x^{-1}A+C]B$. By \cite[Mem. I, p. 220]{L-D}, this has a
solution, valid on $0<\rho_{1} <|x|<\rho_{2}<1$, as a complicated power series
in $x^k\ln^j x$ $(k\in\BbZ, j\in\BbZ^+\cup\{0\})$, with matrix
coefficients. The solution is convergent for \ Re$(t)<0$. This suggests
looking for a solution of (\ref{e:five}) either of the form $\sum
B_{kj}e^{kt}t^j$ or of the form $\sum B_n(e^t-1)^n$. The latter will
give a solution near  $t=0$; recall that in quantum cohomology
solutions are formally constructed in a neighborhood of $t^1
=\ldots=t^k=0$. We can also address convergence issues near a
fixed value of $e^t =\alpha$, by looking for a solution
of the form $\sum B_n(e^t -\alpha )^n$.  In fact, each substitution
gives different information.

\medskip
\noindent {\bf Substitution I:}   We plug  $\sum_{n\ge 0} B_n(e^t-1)^n$ into
(\ref{e:five}), noting that $B(0)=$Id implies $B_0=$Id. We obtain
$$\sum n B_n(e^t -1)^{n-1}e^t =\sum AB_n (e^t -1)^n +\sum CB_n
(e^t-1)^n e^t,
$$
or
\begin{eqnarray*}\lefteqn{\sum n B_n(e^t-1)^n +\sum n B_n(e^t-1)^{n-1}
=}\\
&&\sum AB_n (e^t-1)^n 
+\sum CB_n (e^t-1)^{n+1}
 +\sum CB_n 
(e^t-1)^n.\end{eqnarray*}
This gives the recursion relation
\begin{eqnarray}
B_{n+1} &=&\frac{1}{n+1} [ (A+C-n)B_n + CB_{n-1}],\ n\ge 2,
\label{e:six}\\
B_0 &=& {\rm Id},\  \ B_1=A+C.\nonumber\end{eqnarray}
This two-term relation can be encoded as
\begin{eqnarray*}
\left( \begin{array}{c}
B_{n+1} \\ B_n \end{array}\right)
 &=& 
 \frac{1}{n+1}\left(
\begin{array}{cc} A+C-n& C\\
n+1 & 0
\end{array} \right)
\left( \begin{array}{c}B_{n} \\ B_{n-1}\end{array}\right)  \\
&=& \frac{1}{n(n+1)}\left(
\begin{array}{cc}
A+C-n&C \\
n+1 & 0
\end{array}\right)\\
&&\qquad \cdot \left( \begin{array}{cc}
A+C-(n-1) & C \\
n & 0
\end{array}\right)\left(\begin{array}{c} B_{n-1}\\
B_{n-2}\end{array}\right) \\
&=& \ldots \ = \frac{1}{(n+1)!}\prod_{j=1}^n
\left( \begin{array}{cc}
A+C-j & C \\
j+1 & 0
\end{array}\right)\left(\begin{array}{c} B_1\\B_0\end{array}\right),
\end{eqnarray*}
with the convention that the $j^{\rm th}$ matrix is to the right of the
$(j+1)^{\rm st}$ matrix. This can be rewritten as
\begin{eqnarray*}
B_{n+1} &=&\frac{1}{(n+1)!} \left[ \prod_{j=1}^n\left( \begin{array}{cc}
A+C-j & C \\
j+1 & 0
\end{array}\right) \left( \begin{array}{cc}
A+C & C \\
1 & 0
\end{array}\right)\right]_{(1,1)}
\\
&=& \frac{1}{(n+1)!} \left[ \prod_{j=0}^n \left( \begin{array}{cc}
A+C-j & C \\
j+1 & 0
\end{array}\right)\right]_{(1,1)}.
\end{eqnarray*}

As a check, when $C=0$ we must get $B(t)=e^{tA}$. For fixed $y$, the
Taylor series for $x^y$ at $x=1$ is
$$x^y =1+\sum_{n\ge 1} \frac{1}{n!} \left[ \prod_{j=0}^{n-1}
(y-j)\right] (x-1)^n,
$$
and so
$$e^{tA} ={\rm Id} +\sum_{n\ge 1} \frac{1}{n!}\left[ \prod_{y=0}^{n-1}
(A-j)\right] (e^t-1)^n.
$$
For $C=0,\ B_n=\frac{1}{n!} \prod_{j=0}^{n-1}(A-j)$, so this case
checks. It is also interesting to check the case when $A,C\in\bbC$ are
$1\times 1$ matrices. The explicit solution to (\ref{e:five}) is $B(t)=e^{tA}
e^{(e^t-1)C}$, so we get the identity
$$x^a e^{(x-1)c}=1+\sum_{n\ge 1} \frac{1}{n!}\left[ \prod_{j=0}^{n-1}
\left( \begin{array}{cc}
a+c-j & c\\
j+1 & 0
\end{array}\right)\right]_{(1,1)} (x-1)^n. \\
$$
Writing $x^a=((x-1)+1)^a$ and using the binomial expansion, we get 
$$\sum_{k+\ell =n}  \binom{a}{k} \frac{c^\ell}{\ell
!} =\frac{1}{n!}\left[ \prod_{j=0}^{n-1}\left( \begin{array}{cc}
a+c-j & c\\
j+1 & 0
\end{array}\right)\right]_{(1,1)}.
$$
In summary, Substitution I leads to cute identities.

\bigskip
\noindent {\bf Substitution II:}  For convergence issues, we fix $t_0\in\bbC$,
set $e^{t_0}=\alpha$, and substitute $\sum B_n(e^t-\alpha )^n$ into
(\ref{e:five}). As above, this leads to the recursion relation
$$B_{n+1} =\frac{1}{(n+1)\alpha} [(A+C\alpha -n)B_n + CB_{n-1}],
$$
and so
$$B_{n+1} =\frac{1}{(n+1)! \alpha^{n+1}}\left[ 
\prod_{j=0}^n \left( \begin{array}{cc}
A+C\alpha -j & C \\
(j+1)\alpha & 0
\end{array}\right)\right]_{(1,1)}.$$
Thus
\begin{equation}
B={\rm Id}+\sum_{n\ge 1} \frac{1}{n!} \left[ \prod_{j=0}^{n-1}
\left( 
\begin{array}{cc}
A+e^{t_0}C-j & C   \\
e^{t_0}(j+1) & 0
\end{array}\right)\right]_{(1,1)} (e^{t-t_0}-1)^n.
\label{e:seven}
\end{equation}
The $(1,1)$ entry of the  matrix coefficient of $(e^{t-t_0}-1)^n$
consists of $n$ terms, each
of which is a product of $k$ terms of the form $A+e^{t_0}C-j$, $\ell$
occurrences of $C$, and $m$ terms of the form $e^{t_0}(j+1)$, with
$k+\ell +m\le n$. Thus the supremum norm of the $(1,1)$ entry is
bounded above by $c\cdot n\cdot n!$, where $c$ is a constant depending
on $t_0$ and the sup norms of $A$ and $C_j$, and below by $c^\prime
\cdot n\cdot n!$ for another constant $c^\prime$. The ratio test implies
that the right hand side of (\ref{e:seven}) converges if
$|e^{t-t_0}-1| <1,$
or equivalently if Re$(t-t_0)<2\cos {\rm Im}(t-t_0)$. In particular,
the series converges uniformly on a ball around $t_0$ to a solution of
(\ref{e:five}). 

\bigskip
\noindent 
{\bf Substitution III:}\ \ We substitute $\sum B_{kj}e^{tk}t^j$ into
(\ref{e:five}), with $k\ge 0, j\ge 0, k,j\in\BbZ$. 
The omission of terms with
$k<0$ is motivated by the fact that it works.  From
$B(0)=$ Id, we get $\sum_k B_{k0}=$Id, so assume
$B_{k0}=\delta_{k0}\cdot$Id. The substitution gives
\begin{equation}
\sum k B_{kj}e^{tk}t^j +\sum (j+1)B_{k,j+1}e^{tk}t^j=\sum
AB_{kj}e^{tk}t^j +\sum C B_{k-1,j} e^{tk}t^j,   
\label{e:eight}
\end{equation}
so
\begin{equation}
B_{k,j+1} =\frac{1}{j+1} [(A-k)B_{kj} +C B_{k-1,j}]. 
\label{e:nine}
\end{equation}
For $j>0$, this encodes as
\begin{eqnarray*} \lefteqn{ \left(  \begin{array}{c}
B_{kj}\\  B_{k-1,j}\\ 
\vdots\\ \vdots\\ B_{1j}\\ B_{0j}\end{array}\right) =}\\
&&\\
&& \frac{1}{j}
\left( 
\begin{array}{cccccc}
A-k&C \\
 & A-k+1&C \\
& &\ddots& \ddots\\
&& & \ddots & \ddots \\
&&&  & A-1&C \\
&&&& & A
\end{array}\right) \left( \begin{array}{c}
B_{k,j-1}\\ B_{k-1,j-1}\\ \vdots \\ \vdots\\
B_{1,j-1}\\ B_{0,j-1}) \end{array}\right) \\
&&\\
&=&\ldots 
=  \frac{1}{j!}
\left( 
\begin{array}{cccccc}
A-k&C \\
 & A-k+1&C \\
& &\ddots& \ddots\\
&& & \ddots & \ddots \\
&&&  & A-1&C \\
&&&& & A
\end{array}\right)^j \left( \begin{array}{c}
B_{k0}\\  B_{k-1,0}\\ \vdots \\ \vdots \\
 B_{10}\\ B_{00}\end{array}\right).
\end{eqnarray*}
Since $B_{k0}=\delta_{k0}\cdot$  Id, we get for $k\ge 0,\  j>0$,
\begin{equation}
B_{kj}=\frac{1}{j!} \left[
\left( 
\begin{array}{cccccc}
A-k&C \\
 & A-k+1&C \\
& &\ddots& \ddots\\
&& & \ddots & \ddots \\
&&&  & A-1&C \\
&&&& & A
\end{array}\right)
 ^j\right]_{(1,k+1)}.
\label{e:ten}
\end{equation}
Note that the assumption $B_{kj}=0$ for $k<0$ is consistent with
(\ref{e:nine}). The theory in \cite{L-D} gives uniform convergence of $\sum
B_{kj}e^{kt}t^j$ to a solution for $-\frac{1}{\epsilon}<$
Re$(t)<-\epsilon <0$ for any $\epsilon >0$, and so it must coincide
with the solution given in Substitution II.  This again leads to
 identities.

(\ref{e:ten}) has an infinite dimensional interpretation. Let
$H=L_{+}^2(S^1, \bbC^m)$ denote the Hilbert space of
$L^2\ \bbC^m$-valued functions on $S^1$ with Fourier expansions
containing only $e^{in\theta}$ with $n\ge 0$. Let $D: H\rightarrow H$
be the first order operator
\begin{equation}\label{d}
D=-\sqrt{-1} \frac{d}{d\theta} +A+e^{-\sqrt{-1}\theta}C.
\end{equation}
(The range of $D$ is $H$ plus the span of $e^{-\sqrt{-1}\theta}$, so
we actually compose $D$ with the projection $L^2\to H.$)
Let $e^{\sqrt{-1}n\theta}_{(j)}$ denote $(0,\ldots
0,e^{\sqrt{-1}n\theta},0,\ldots,0)$, where $e^{\sqrt{-1} n\theta}$
occurs in the $j^{th}$ slot. Then
$De^{\sqrt{-1} n\theta}_{(j)}
=ne_{(j)}^{\sqrt{-1}n\theta}+A_{jk}e_{(k)}^{\sqrt{-1}
n\theta}+C_{jk}e_{(k)}^{\sqrt{-1}(n+1)\theta}$, so in the basis
\hfill\break 
$\{
e_{(j)}^{\sqrt{-1} n\theta}\}_{n\ge 0,j=1,\ldots,m}$, $D$ has matrix
$$\left( 
\begin{array}{ccccc}
A&C \\
 & A+1&C \\
&  & A+2&C\\
&&& \ddots &\ddots
\end{array}\right)
$$
which we also denote by $D$.

\medskip
We compute
\begin{eqnarray*}
B &=& \sum_{k,j\ge 0} \frac{1}{j!} \left[ 
\left( \begin{array}{cccccc}
A-k&C \\
 & A-k+1&C \\
& &\ddots& \ddots\\
&& & \ddots & \ddots \\
&&&  & A-1&C \\
&&&& & A
\end{array}\right)
 ^j\right]_{(1,k+1)}\cdot e^{kt}t^j \\
&=&\sum_{k\ge 0} \left[ \exp\left[ t\left( 
\begin{array}{cccccc}
A-k&C \\
 & A-k+1&C \\
& &\ddots& \ddots\\
&& & \ddots & \ddots \\
&&&  & A-1&C \\
&&&& & A
\end{array}\right)
\right] \right]_{(1,k+1)}e^{kt} \\
&=& \sum_{k\ge 0} \left[ \exp\left[ t\left( 
\begin{array}{cccccc}
A-k&C \\
 & A-k+1&C \\
& &\ddots& \ddots\\
&& & \ddots & \ddots \\
&&&  & A-1&C \\
&&&& & A
\end{array}\right)
\right]\right]_{(1,k+1)}
\end{eqnarray*}

It is easy to check that the $(1,k+1)$ entry of this last matrix,
denoted $\exp [t D_k]$, equals the $(1,k+1)$ entry of $\exp [tD]$ by
comparing entries for $D_k^j$ and $D^j$. This gives:
\begin{proposition} The solution of $\dot B(t) = (A + e^tC)B(t), B(0)
= {\rm Id}$ is
given by
\begin{equation}
B=\sum_{k\ge 0} [\exp (tD)]_{(1,k+1)}. 
\label{e:eleven}
\end{equation}
with $D$ given by (\ref{d}).\end{proposition}

In particular,
$$
B_{ij} =\sum_{k\ge 0} \langle \exp (tD)(e_{(j)}^{\sqrt{-1}k\theta}),\
1_{(i)} \rangle,  $$
since $e^{\sqrt{-1}0\theta}=1$. 
We  can
extend $D$ by its matrix representation to act on \\
$L_{+}^2$ ($S^1,
M_{m\times m}(\bbC)$), and then
\begin{equation}
B=\sum_{k\ge 0} \langle \exp (tD)(e^{\sqrt{-1}k\theta}{\rm Id}),\ 
{\rm Id}\rangle .  
\label{e:thirteen}
\end{equation}
As in \S2, we can now continue to solve in other $H^2$ directions.

For ${\mathbb P}^m$, $B$ can be computed from (\ref{e:eleven}) since
$A,C$ are particularly simple.
For other spaces,
we can make (\ref{e:thirteen})
 more explicit by attempting to diagonalize $D$. Let
$$e_k =\left( 
\begin{array}{cccccccc}
0&\cdots &0&
{\rm Id} &-C & C^2/2 &
-C^3/3! &\cdots \end{array}\right),$$
with  Id  in the $k^{\rm th}$ slot. The $\{ e_k\}$ are linearly
independent, and
 $ e_kD =e_k (A+k)$.  For
$$E=\left( 
\begin{array}{ccccc}
{\rm Id} &-C&C^2/2&C^3/3!&\cdots\\
 & {\rm Id}&-C &C^2/2&\cdots\\
 & & {\rm Id}&-C&\ddots \\
 &&  &\ddots&\ddots \\
 &  &  &&\ddots
\end{array}\right)
$$
we have
$$
ED =\left(
\begin{array}{ccccc}
A&-CA&\frac{C^2A}{2}& -\frac{C^3A}{3!}&\cdots \\
& A+1&-C(A+1) & \frac{C^2}{2}(A+1)&\cdots       \\
 &  & A+2&-C(A+2)&\ddots \\
&  &  & A+3&\ddots \\
&&&& \ddots
\end{array}\right)
$$
\begin{eqnarray*}
&=&\left( \begin{array}{ccccc}
A&-AC -[C,A]& A\frac{C^2}{2}+[\frac{C^2}{2},A]&
-A\frac{C^3}{3!}-[\frac{C^3}{3!},A]&\cdots \\
 & A+1& (A+1)C-[C,A]&(A+1)\frac{C^2}{2}+[A,\frac{C^2}{2}]&\cdots  \\
&&A+2&(A+2)C - [C,A]&\ddots \\
& &&\ddots
&\ddots \\
\end{array}\right) \\
&=& GE+T
\end{eqnarray*}
for
\begin{eqnarray*}G &=& \left( 
 \begin{array}{cccc}
A \\
&A+1 \\
&&A+2 \\
&&&\ddots 
\end{array}
\right),\\
 T &=& \left( 
\begin{array}{ccccc}
 0&-[C,A]&[C^2/2,A]&-[C^3/3,A]&\cdots  \\
 & 0& -[C,A]& [C^2/2,A]&\cdots \\
 &  & 0&-[C,A]& \ddots \\
&&&\ddots&\ddots \\
\end{array}\right).
\end{eqnarray*}
Since $E$ is invertible,
\begin{equation}
D=E^{-1}(G+TE^{-1})E, 
\label{e:fourteen}
\end{equation}
so $TE^{-1}$ measures the obstruction to diagonalizing $D$ due to the
noncommuting of $A$ and $C$.  Note that
$$E^{-1}= \left(
 \begin{array}{ccccc}
{\rm Id}&C&C^2/2&C^3/3!&\cdots \\
& {\rm Id}&C&C^2/2&\cdots \\
&&{\rm Id}&C&\ddots \\
&&&\ddots&\ddots \\
\end{array}\right).
$$
By (\ref{e:thirteen}), (\ref{e:fourteen}),
\begin{eqnarray*}
B &=&\sum_{k\geq 0} \langle E^{-1}\ e^{t(G+TE^{-1})}E
e^{\sqrt{-1}k\theta}{\rm Id},{\rm Id}\rangle\\
  &=&\sum_{k\geq 0}\left\langle
   (E^{-1})^T{\rm Id},
   e^{t(G+TE^{-1})} E
e^{\sqrt{-1}k\theta}{\rm Id}
  \right\rangle {}^{{}^{ {}^ { {}^
  {{}^{{}^{\overline{\ \ }}}}}}}
\end{eqnarray*}
where the final bar denotes complex conjugate.  This uses the identity
$\langle Av,w\rangle = \overline{\langle A^Tw,v\rangle}$ for a real
endomorphism $A$ of a complex noncommutative algebra, provided the
components of $v,w$ satisfy $[v_i,\bar w_i]=0.$
Since $Ee^{\sqrt{-1}k\theta}$ Id  is the k${}^{\rm th}$ column of
$E$, we get
$$\overline B=\sum_{k\geq 0} \Bigg\langle 
 \left(
 \begin{array}{c}
{\rm Id} \\
C \\
C^2/2 \\
\vdots
\end{array}\right),
 e^{t(G+TE^{-1})}
\left( 
\begin{array}{c}
(-1)^{k-1} C^{k-1}/(k-1)! \\
(-1)^{k-2}  C^{k-2}/(k-2)! \\
\vdots \\
{\rm Id} \\
0 \\
\vdots
\end{array}\right)\Bigg\rangle.
$$
Thus
\begin{eqnarray} \label{e:fifteen}
\overline B &=& \Bigg\langle \left( 
 \begin{array}{c}
{\rm Id} \\
C \\
C^2/2 \\
\vdots
\end{array}\right),
e^{t(G+TE^{-1})}\left(
\begin{array}{c}
e^{-C} \\
e^{-C} \\
e^{-C} \\
\vdots
\end{array}\right) \
 \Bigg\rangle.  \\
&=& \left\langle e^{Ce^{\sqrt{-1}\theta}}, e^{t(G+TE^{-1})}\sum_{k\geq
0}e^{-C} e^{\sqrt{-1}k\theta}\right\rangle.\nonumber
\end{eqnarray}
This expression will be more useful than (\ref{e:eleven}) if $[A,C^k]$
and hence $TE^{-1}$
is sparse.

\section{Flat sections for $\BbP^m$}

In computing  the flat sections for
$\BbP^m$, the only  step is the ODE for $H^2$. Here $A_j^i$
equals 
\begin{eqnarray*}
 &&
\left( \int_{\BbP^n} T_1\cup T_j\cup T_s\right) h^{si}=\left[ 
\left( 
\begin{array}{cccc}
&&1&0 \\
&\cdot^{\ \cdot^{\ \cdot}}&\cdot^{\ \cdot^{\ \cdot}} \\
1&\cdot^{\ \cdot^{\ \cdot}} \\
0
\end{array}\right)
\left( \begin{array}{cccc}
& & &1 \\
&& \cdot^{\ \cdot^{\ \cdot}}& \\
& \cdot^{\ \cdot^{\ \cdot}}&& \\
1&&&
\end{array}\right)\right]^i_j \\
&=&
\left( 
\begin{array}{cccc}
0&1 \\
&\ddots&\ddots \\
&&\ddots&1 \\
&&&0 
\end{array}\right)_j^i
\end{eqnarray*}
and
$$C=\left( 
\begin{array}{ccc}
&& \\
&& \\
&&1
\end{array}\right)
 \left(
\begin{array}{ccc}
&&1 \\
&\cdot^{\ \cdot^{\ \cdot}} \\
1
\end{array}\right)
=\left( 
\begin{array}{ccc}
&&\\
&& \\
 1&&
\end{array}\right), 
$$
where all empty slots are $0$. Recall that we solved (\ref{e:five})
for $-A, -C$.

As in (\ref{e:eleven}), we have
\begin{equation}\label{four.one} B = \sum_{k,n\geq 0}\frac{t^n}{n!}
D^n_{(1,k+1)}.\end{equation}
Note that the $k=0$ summand is 
$$\sum_n \frac{t^n}{n!}(-A)^n
= \left(\begin{array}{ccccc}
1&-t&t^2/2&\cdots&(-1)^mt^m/m!\\
&1&-t&\ddots&\vdots\\
&&\ddots&\ddots&t^2/2\\
&&&\ddots&-t\\
&&&&1\end{array}\right),$$
and that the only other contribution from the sum for $n=0$ or $n=1$
is $tC$ when $(n,k) = (1,1).$ Thus we may assume below that $n>1$ and
$k>0.$  

Write $D = \ma+\mc$, where 
$$\ma = \left(\begin{array}{cccc}
-A&&&\\
&-A+1&&\\
&&-A+2&\\
&&&\ddots\end{array}\right),
\ 
\mc = \left(\begin{array}{cccc}
0&-C&&    \\
&0& -C&   \\
&&0&\ddots   \\
&&&\ddots  \end{array}\right).
$$
Let $P(a,b)=\{(r_1,\ldots,r_a):r_i\in \BbZ^+, \sum r_i =b\}$ be the
set of unordered partitions of $b$ into $a$ terms. 
Since $C^2=0$ and hence $\mc^2=0$, we have 
\begin{eqnarray} \label{four.two}
D^n &=&
\sum_s\sum_{P(s,n-s)}\ma^{r_1}\mc\ma^{r_2}\mc\cdot\ldots\cdot
\ma^{r_s}\mc\nonumber\\ 
&&\qquad + \sum_s\sum_{P(s,n-s+1)}\ma^{r_1}\mc\ma^{r_2}\mc\cdot\ldots\cdot
\ma^{r_s}\nonumber\\  
&&\qquad + \sum_s\sum_{P(s,n-s)}\mc\ma^{r_1}\mc\ma^{r_2}\mc\cdot\ldots\cdot
\ma^{r_s}\\  
&&\qquad + \sum_s\sum_{P(s,n-s-1)}\mc\ma^{r_1}\mc\ma^{r_2}\mc\cdot\ldots\cdot
\ma^{r_s}\mc\nonumber\\  
&=& \rone + \rii +\riii + \riv.\nonumber\end{eqnarray}
 A fixed partition contributes only
 one nonzero entry in the first
row of I${}^{(n)}$, namely 
$(-1)^s(-A)^{r_1}C(-A+1)^{r_2}C\cdot\ldots\cdot(-A+s-1)^{r_s}C$
in the $(1,s+1)$ slot.  Thus
$$
\rone_{(1,k+1)} =
\sum_{P(k,n-k)}(-1)^k(-A)^{r_1}C(-A+1)^{r_2}C\cdot\ldots\cdot(-A+k-1)
^{r_k}C.$$ 
Since $(-A+c)^r = \sum_{u=0}^r(-1)^u \binom{r}{u}A^uc^{r-u} = \sum_{u=0}^r
(-1)^u\binom{r}{u}(\delta_{i,j-u})c^{r-u} $, and since the product of
$n\times n$ matrices
$X^{(1)}CX^{(2)}C\cdot\ldots\cdot X^{(k)}C$ has only nonzero entries
in the first column, with $(\ell,1)$ entry $x_{\ell
n}^{(1)}x_{1n}^{(2)}x_{1n}^{(3)}\cdot\ldots\cdot x_{1n}^{(k)}$, we get
\begin{eqnarray}\label{four.three}
\rone_{(1,k+1)} &=& \sum_{P(k,n-k)} 
(-1)^{(k-1)m+k}(-1)^{r_1}\binom{r_2}{m}1^{r_2-m}\binom{r_3}{m}2^{r_3-m}
\cdot\ldots\\
&&\qquad\nonumber  \cdot
\binom{r_{k}}{m} (k-1)^{r_{k}-m}E^{m+1-r_1}_1,
\end{eqnarray}
where $E^i_j$ is the $(m+1)\times (m+1)$ matrix with one in the
$(i,j)$ entry and all other entries zero.
 Here we set $\binom{r}{m}=0$ if $r<m.$  Notice that $E^{m+1}_1$ does not
occur in (\ref{four.three}).

$\riv$ has a nonzero entry in the $(1,k+1)$ slot iff
$\ma^{r_1}\mc\cdot\ldots\cdot\ma^{r_s}\mc$ has a nonzero entry in the
$(2,k)$ slot.  This entry is $(-1)^{k}C(-A+1)^{r_1}C\cdot\ldots\cdot
(-A+k-1)^{r_{k-1}} C$, which has a nonzero entry only in the $(1,m+1)$
slot. As above, we get
\begin{eqnarray}\label{four.three_a}
\riv_{(1,k+1)} &=& \sum_{P(k-1,n-k)}
(-1)^{km+k-m}\binom{r_1}{m}1^{r_1-m}
\binom{r_2}{m}2^{r_2-m}\\
&&\qquad\nonumber\binom{r_3}{m}3^{r_3-m}
\cdot\ldots \cdot
\binom{r_{k-1}}{m} (k-1)^{r_{k-1}-m}E^{m+1}_1.\end{eqnarray}
Letting $r_1+1=j$ run over $\{1,\ldots,m+1\}$
and then relabeling $(r_2,\ldots,r_{k})$ as $(r_1,\ldots,r_{k-1})$ in
(\ref{four.three}), we can combine
(\ref{four.three}), (\ref{four.three_a}) to get
\begin{eqnarray}\label{four.three_b} 
\rone_{(1,k+1)} + \riv_{(1,k+1)} &=& \sum_{j=1}^{m+1} \sum_{P(k-1,
n-k-j+1)} (-1)^{(k-1)m+k+j-1} \binom{r_1}{m}1^{r_1-m}\nonumber \\
&&\qquad \cdot\binom{r_2}{m}2^{r_2-m}
\cdot\ldots  \cdot
\binom{r_{k-1}}{m} (k-1)^{r_{k-1}-m}E^{m+2-j}_1.\end{eqnarray}

The other terms are handled similarly.  
We have
$$
\rii_{(1,k+1)} = \sum_{P(k+1,n-k)}
(-1)^k(-A)^{r_1}C(-A+1)^{r_2}C\cdot\ldots\cdot
C(-A+k)^{r_{k+1}}.                     
$$
The $(i,\ell)$ entry of $X^{(1)}CX^{(2)}C\cdot\ldots\cdot
X^{(k)}CX^{(k+1)}$ is
$x_{in}^{(1)}x_{1n}^{(2)}x_{1n}^{(3)}\cdot\ldots\cdot
x_{1n}^{(k)}x_{1\ell}^{(k+1)}.$  Since $((-A)^{r_{1}})^i_{m+1} =
(-1)^{r_1}\delta^i_{m+1-r_1}$,  we get
\begin{eqnarray}\label{four.five}
\rii_{(1,k+1)} &=& \sum_{P(k+1,n-k)}
\sum_{\ell =1}^{m+1} (-1)^{(k-1)m+\ell+k+r_1-1} 
\binom{r_2}{m}1^{r_2-m}\nonumber \\
&&\qquad \cdot\binom{r_3}{m}2^{r_3-m}
\cdot\ldots  \cdot
\binom{r_{k}}{m} (k-1)^{r_{k}-m}\binom{r_{k+1}}{\ell-1}k^{r_{k+1}-\ell+1}
E^{m-r_1+1}_\ell.\end{eqnarray}
Since $r_1>0$, (\ref{four.five}) has no entries in the $m+1^{\rm st}$ row.
As with $\riv$, this missing
row appears in $\riii$.  Setting  $r_1+1=j$ and shifting indices on
the $r_i$ as above,
we get
\begin{eqnarray} \label{four.six}
\rii_{(1,k+1)} + \riii_{(1,k+1)} &=& \sum_{j=1}^{m+1} \sum_{P(k,
n-k-j+1)}\sum_{\ell =1}^{m+1} (-1)^{(k-1)m + \ell+k +j-2}
\binom{r_1}{m}1^{r_1-m}\nonumber \\
&&\qquad \cdot\binom{r_2}{m}2^{r_2-m}
\cdot\ldots  \cdot
\binom{r_{k-1}}{m} (k-1)^{r_{k-1}-m}\binom{r_{k}}{\ell-1}\\
&&\qquad \cdot k^{r_{k}-\ell+1}
E^{m-j+2}_{\ell}.\nonumber
\end{eqnarray}

Finally, finding the flat sections for $\nabla_\hbar$ as in (3.2)
requires replacing $A,C$ with
$(-\sqrt{-1}\hbar)^{-1}A,(-\sqrt{-1}\hbar)^{-1}C.$  This multiplies
the entry $(-1)^u\binom{r}{u}c^{r-u}$ in $(-A+c)^r$ by $\scl^{-u}$.  Since
there are $k-1$ terms
with binomial coefficient $\binom{r_i}{m}$, one term of the form
$(-A)^{j-1}$, and 
$k$ copies of $C$ in each term on the right hand side of 
(\ref{four.three_b}), this term
 is multiplied by
$\scl^{-((k-1)m+k+j-1)}.$
 Similarly, each entry in (\ref{four.six}) and
 is multiplied by $\scl^{-(mk+k+j-1)}.$  

Set $\alpha = \scl.$
Combining (\ref{four.one}), (\ref{four.two}), (\ref{four.three_b}),
(\ref{four.six}) (and remembering the
contribution from $tC$), we obtain:
\begin{theorem}\label{below}
  A basis of the flat sections for the Dubrovin
connection 
$\nabla_\hbar= \nabla_{\sqrt{-1}\alpha}$ on
$T{\mathcal H}|_{H^2({\mathbb P}^m)},$ is given by the columns of
\begin{eqnarray}\label{pnsections}B(t,\hbar) &=& A(t,\hbar)+
\sum_{n\geq 2, k\geq 1}\frac{t^n}{n!}\sum_{j=1}^{m+1}\left[
\alpha^{-((k-1)m
+k+j-1)}\right. \nonumber\\
&&\qquad\ 
\cdot  \sum_{P(k-1,
n-k-j+1)} (-1)^{(k-1)m+k+j-1}
\left(\prod_{s=1}^{k-1}\binom{r_s}{m}s^{r_s-m}\right) E_1^{m+2-j}\nonumber \\
&&\qquad    
+\alpha^{-(mk+k+j-1)} \\
&&\qquad\  \cdot
\sum_{P(k,
n-k-j+1)}\sum_{\ell =1}^{m+1} (-1)^{(k-1)m+\ell+k+j-2}
\left(\prod_{s=1}^{k-1}\binom{r_s}{m}s^{r_s-m}\right)\nonumber\\
&&\left. \qquad\  \cdot
\binom{r_k}{\ell-1}k^{r_k-\ell+1}E_{\ell}^{m-j+2}\right],\nonumber\end{eqnarray}
with 
\begin{eqnarray*}\lefteqn{A(t,\hbar ) =} \\
&&\!\!\!\!\!\!\!\!
\left(\begin{array}{ccccc}
1&\alpha^{-1}t&\alpha^{-2}t^2/2&\cdots&\alpha^{-m}t^m/m!\\
0&1&\alpha^{-1}t&\ddots&\vdots\\
\vdots&\ddots&\ddots&\ddots&\alpha^{-2}t^2/2\\
0&&\ddots&\ddots&\alpha^{-1}t\\
\alpha^{-1}t&0&\cdots&0&1    \end{array}\right).        \end{eqnarray*}
\end{theorem}

\begin{remark}  
In non-$H^2$ directions,
there are no quantum
corrections to the cup product; i.e.~the quantum product is
$t$-independent in these directions. Specifically,
for $T_\ell \not= T_1$, we have $C=0$ and the corresponding $A$ matrix is
\begin{eqnarray*}
A_j^i = \left( \int_{\BbP^n} T_\ell \cup T_j\cup T_s\right) h^{si} &=&
(\delta_{j+s}^{n-\ell}) (\delta^k_{s}\delta_{n-k}^{i}) \\
&=&\delta_{\ell +j}^i.
\end{eqnarray*}
Using Proposition 2.2, we see that the flat sections for the Dubrovin
connection over all of $\mathcal H$ are given by
$$   \exp \left[ -\sqrt{-1}\sum^m_
  {\scriptstyle{{{\ell =0},
  {{\ell\not= 1}}}}} 
  t^\ell \left(\frac{\sqrt{-1}}{\hbar}\right)\Gamma_\ell\right]
B(t^1),$$
with
  $(\Gamma_\ell )_j^i =\delta_{\ell +j}^i$ and $B(t)$ 
the matrix in (\ref{pnsections}).  The lack of
significant information in non-$H^2$ directions is 
a special feature of ${\mathbb P}^m.$
\end{remark}

It is interesting to compare   (\ref{pnsections}) to the basis of
 flat sections obtained in
\cite{P}.   As in \cite{P}, for
 a vector field $f^iT_i$ in $i^*T{\mathcal H}$, 
the equation (\ref{dub}) is equivalent to the
 system 
\begin{eqnarray}\label{pansys}
-\hbar^{-1} \partial_t(f^i) &=& f^{i-1}, \ i>0,\\
-\hbar^{-1} \partial_t(f^0) &=& e^{t}f^m,\nonumber \end{eqnarray}
where $t = t^1$ and 
$\partial_t$ denotes differentiation by the vector field $T_1.$
(Our choice of $\hbar$ equals $-\hbar^{-1}$ in \cite{P}).  Thus $f^m$
determines the other $f^i$, and $f^m$ must satisfy
\begin{equation}\label{deqn}
(-\hbar^{-1}\partial_t)^{m+1}f^m - e^{t}f^m =0.\end{equation}
Letting $H$ be a formal variable (which we think of as
the hyperplane class in $H^2({\mathbb P}^n)$),
it is easy to check that 
$$S = \sum_{d\geq 0}\frac{e^{(-\hbar H + d)t}}{\prod_{r=1}^d
(H-\hbar^{-1}r)^{m+1}} \ \ {\rm mod} \ H^{m+1}$$
formally solves (\ref{deqn}). (The denominator is one for $d=0.$)
 For $S = \sum_{b=0}^m S_bH^{m-b}$, each
$S_b = (\partial^{m-b}/\partial H^{m-b})|_{H=0} S$ must also solve
(\ref{deqn}).  For a fixed $S_b$, a solution of (\ref{pansys}) is then 
$\partial_t^{m-s}S_b T_s.$  In other words, a basis of the flat sections
is given by the columns of the matrix $M=M(t,\hbar)$ with
$$M_b^s = (-\hbar^{-1}\partial_t)^{m-s}(\partial^{m-b}/\partial
H^{m-b})|_{H=0} S.$$
The uniqueness of flat sections with initial condition gives
\begin{equation} M(0,\hbar)B(t,\hbar)
=M(t,\hbar).
\label{e:twenty-nine}
\end{equation}
This produces complicated identities in $t,\hbar.$

Another approach for ${\mathbb P}^1$ is in \cite{V}.


\section{Remarks on the general case}

In the general case of the small quantum product, the
analogue of the nontrivial ODE (\ref{e:five}) is of the form
\begin{equation}
\dot B(t)=B(t) (A+\sum_\beta \exp[\sum_{i=1}^p t^i
a_{i\beta}]C_\beta), 
\label{e:thirty}
\end{equation}
where $a_{i\beta}=\int_\beta T_i$.  We solve this by Substitution
III. Namely, we first set \hfill\break
$(t^1, \ldots, t^p)=(t,0,\ldots,0)$ as
in  (\ref{e:five}) and assume $B(t)=\sum B_{kj} e^{kt}t^j$.
(\ref{e:eight}) becomes
\begin{equation}
\sum kB_{kj} e^{tk}t^j +\sum (j+1)B_{k,j+1} e^{tk}t^j=\sum
B_{kj}Ae^{tk} t^j +\sum_{i=1}^n B_{k-a_i} C_i e^{kt}t^j,
\label{e:thirty-one}
\end{equation}
for $a_i=a_{1\beta_i},\ c_i=c_{\beta_i}$, where $\{ \beta_i\}$ is the
set of $\beta$'s with $I_\beta(T_j,T_1,T_\ell )\not= 0$ for some
$j,\ell$.  (\ref{e:nine})
 becomes
$$B_{k,j+1} =\frac{1}{j+1} [B_{kj}(A-k)+\sum_{i=1}^n B_{k-a_i,j}C_i]
$$
If $M$ is weakly monotone,  $a_i>0$, and
 the expression for $j!B_{kj}$ in (the transposed version of)
(\ref{e:ten}) becomes
$$
\left[ \left( 
\begin{array}{ccccccccc}
A-k \\
\vdots&A-k+1 \\
C_1&\vdots&A-k+2 \\ 
\vdots&C_1&\vdots&\ddots \\
C_2&\vdots&C_1&&\ddots \\
\vdots&C_2&\vdots&\ddots&&\ddots \\
C_n&\vdots&C_2&&\ddots&&\ddots \\
\vdots&C_n&\vdots&\ddots&&\ddots&&A-1 \\
\vdots&\vdots&C_n&&\ddots&&\ddots&(C_1)&A
\end{array}\right)^j\right]_{(k+1,1)}
$$
where $(C_1)=C_1$ if $a_1=1$ and 0 otherwise. In all columns, the
number of zeros between the $A-k+...$ and $C_1$ is $a_1$, and the
number of zeros between $C_i$ and $C_{i+1}$ is $a_i$.
Setting
$$D=-\sqrt{-1}\frac{d}{d\theta}+A+\sum_{i=1}^n
e^{\sqrt{-1}a_i\theta}C_i,
$$
we get
$$B(t) =\sum_{k\ge 0} \langle \exp (tD){\rm Id}, e^{\sqrt{-1}k\theta}
{\rm Id}\rangle .
$$
Here $A_j^i =(\int_M T_1\cup T_j\cup T_\ell )h^{\ell i}=\langle
T_1\cup T_j,T_\ell\rangle h^{\ell i}=(T_1\cup T_j)^i$ where $T_1\cup
T_j=(T_1\cup T_j)^i T_i$. For the next step in the Frobenius
Integrability Theorem, we solve
\begin{eqnarray*}\dot B_{(2)}(t) &=&
B_{(2)}(t) \left[ A_{(2)} +\sum_\beta e^{t_0^1\int_\beta
T_1+t\int_\beta T_2} I_\beta (T_j,T_2,T_\ell )h^{\ell r}\right]\\
&=& B_{(2)}(t) [A_{(2)}+\sum_{i=1}^{n_2} e^{a_i^{(2)}t}
C_i^{(2)}],\\
B_{(2)}(0) &=& B(t_0^1),\end{eqnarray*}
with $(C_i^{(2)})_j^r =\exp[t_0^1 \int_{\beta_i} T_1] I_{\beta_i}
(T_j,T_2,T_\ell )h^{\ell_r}$, and $\{ \beta_1,\ldots,\beta_{n_2}\}$ is
the set of $\beta\in H_2$ with $I_\beta (T_j,T_2,T_\ell )\not= 0$ for
some $j,\ell$. Set $D^{(2)}$\hfill\break 
$=-\sqrt{-1}
\frac{d}{d\theta}+A_{(2)}+\sum_{i=1}^{n_2} e^{\sqrt{-1}
a_i^{(2)}\theta}C_i^{(2)}$ with $(A^{(2)})_j^i =(T_2\cup T_j)^i$. Calling
$B=B_{(1)}$ and proceeding, we get that a basis of the flat sections
are the columns of
\begin{equation}
g=
\prod_{i=1}^k \sum_{j\ge 0} \langle \exp (t^iD_{(i)}){\rm Id},
e^{\sqrt{-1} j\theta} {\rm Id}\rangle 
\label{e:thirty-two}
\end{equation}
with
\begin{eqnarray*}
D_{(i)} &=& -\sqrt{-1} \frac{d}{d\theta} +A^{(i)}+\sum_{k=1}^{n_i}
e^{\sqrt{-1} a_k^{(i)}\theta} C_k^{(i)}, \\
(A^{(i)})_j^r &=& (T_i\cup T_j)^r,\ \ \ {\rm for}\ \ T_i\cup T_j
=(T_i\cup T_j)^r T_r, \\
(C_k^{(i)})_j^r &=& \exp \left[ \sum_{\alpha =1}^i t^\alpha
\int_{\beta_k} T_\alpha \right]I_{\beta_k}(T_j,T_i,T_\ell )h^{\ell r},
\end{eqnarray*}
for $\{ \beta_k\}$ the set of $\beta\in H_2$ with
$I_{\beta_k}(T_j,T_i,T_\ell )\not= 0$ for some $i,\ell$. 

We can recover all GW invariants from the solution
(\ref{e:thirty-two}.)
For example, to recover the $C_k^{(1)}$, set $t^1 = t, t^2=\ldots = t^k=0.$
For $C=C_1^{(1)}$, note that
(\ref{e:eight}) implies $B_{1,1}=C$. Since
$$Z(t)\equiv \sum_{k\ge 0} \langle \exp (t D_{(1)}){\rm Id},\
e^{\sqrt{-1}k\theta} {\rm Id}\rangle =\sum B_{k,j}e^{kt} t^j,
$$
we get
$$C_1^{(1)} =(t^{-1}e^{-t}(Z(t)-{\rm Id}))|_{t=0}
=(t^{-1}e^{-t}(g(t,0,\ldots,0)-{\rm Id}))|_{t=0}.
$$
We then multiply $C_1^{(1)}$ by $h_{ij}\cdot\exp
[-t\int_{\beta_1}T_1]$ to recover $I_{\beta_1}(T_j,T_1,T_\ell )$. We
similarly recover $C_k^{(1)}$, and then the $C_k^{(i)}$.

We conclude with some remarks on the big quantum product.
As in \cite{BCPP}, the big quantum product 
is an associative, potential product
and hence produces a flat connection.  To find flat sections, 
we must now solve the formal ODE
$$\frac{\partial g_j^i}{\partial t} =\sqrt{-1}g_r^i \left[ T_j\cdot
T_1\cdot T_\ell +\sum_n \frac{t^n}{n!}\sum_\beta I_\beta
(T_j,T_1,T_\ell ,\overbrace{T_0,\ldots,T_0}^n)\right] h^{\ell r}
$$
at the first step. 
The sum over $n,\beta$ is not finite in general. Since there
are no exponential terms, we rewrite this as
$$\dot B(t)=B(t)[A+\sum_n t^n C_n]
$$
and plug in $B(t)=\sum_{k\ge 0}B_k t^k$. The recurrence relation is
$$B_k =\frac{1}{k} \left[ B_{k-1} A+\sum_{n,k} B_{k-n-1} C_n\right]
$$
with $B_0=$ Id,  $B_{-j}=0$ for $j>0$. This encodes as
\begin{eqnarray*}
\left( 
\begin{array}{c}
B_k \\
 B_{k-1} \\
\vdots \\
B_2 \\
B_1 \\
{\rm Id}
\end{array}\right)& =&
\left(
\begin{array}{cccccc}
\frac{A}{k}&\frac{C_1}{k}&\frac{C_2}{k}&\ldots&\frac{C_{k-1}}{k}&
0 \\
1&0&0&\cdots&0&0 \\
0&1&0&\cdots&0&0 \\
\vdots&&\ddots&&\vdots&\vdots \\
0&0&\cdots&1&0&0 \\
0&0&0&\cdots&1&0
\end{array}\right) \left( 
\begin{array}{c}
B_{k-1} \\
B_{k-2} \\
\vdots \\
B_1 \\
{\rm Id} \\
0
\end{array}\right) \\
&&\\
&=&\ldots\\
&&\\
&=&\frac{1}{k!}
\left( 
\begin{array}{cccccc}
A&C_1&C_2&\cdots&C_{k-1}&0 \\
1&0&0&\cdots&0&0 \\
0&1&0&\cdots&0&0 \\
\vdots&&\ddots&&\vdots&\vdots \\
0&0&\cdots&1&0&0 \\
0&0&0&\cdots&1&0
\end{array}\right)^k
\left( 
\begin{array}{c}
1\\
0 \\
\vdots\cr 0 \\
0 \\
0
\end{array}\right) \\  
&\equiv&\frac{1}{k!} \left( \alpha_{(k)}^k\right)_1^1,
\end{eqnarray*}
As before, setting
$$\alpha =\left( 
\begin{array}{ccccccc}
A&C_1&\cdots&C_{k-1}&C_k&C_{k+1}&\ldots \\
1&0&\ldots&0&0&0&\ldots \\
0&1&\ldots&0&0&0&\ldots \\
0&0&1&0&0&0&\ldots \\
0&0&0&\ddots&0&0&\ldots \\
\vdots&\vdots&\vdots&\vdots&\vdots&\vdots&\cdots
\end{array} \right)
$$
it is easily checked that $\left( \alpha_{(k)}^k\right)_1^1
=(\alpha^k)_1^1$, and so
$$B(t)=\sum B_kt^k =\sum \frac{1}{k!}(\alpha^k)_1^1 t^k
=(e^{t\alpha})_1^1.
$$
We now proceed as before to generate solutions to the big quantum
product by diagonalizing $\alpha$ as much as possible.
These computations, as well as computations for
quantum cohomology coupled to gravity
in genus zero \cite{LT}
will be discussed in future work.

\bibliographystyle{amsplain}

\end{document}